\documentclass[12pt,twoside,leqno]{article}
\usepackage{amsmath}
\usepackage{amssymb}
\usepackage{amsxtra}
\usepackage{amscd}
\usepackage{amsthm}
\usepackage[mathscr]{eucal}
\usepackage{xr}

\swapnumbers

\theoremstyle{plain}
\newtheorem{thm}[subsection]{Theorem}
\newtheorem{prop}[subsection]{Proposition}
\newtheorem{cor}[subsection]{Corollary}
\newtheorem{lem}[subsection]{Lemma}

\theoremstyle{definition}
\newtheorem{defn}[subsection]{Definition}

\newtheorem{rem}[subsection]{Remark}
\newtheorem{para}[subsection]{}

\newtheorem{sbrem}[subsubsection]{Remark}

\newenvironment{pf}{\proof[\proofname]}{\endproof}
\newenvironment{pf*}[1]{\proof[#1]}{\endproof}

\newcommand{\N}{{\mathbb{N}}}
\newcommand{\Q}{{\mathbb{Q}}}
\newcommand{\Z}{{\mathbb{Z}}}

\newcommand{\isom}{\cong}
\newcommand{\Spec}{\operatorname{Spec}}

\newcommand{\fs}{{\rm fs \ }}
\newcommand{\gp}{\mathrm{gp}}
\newcommand{\Gmlog}{\mathbb{G}_{m,\log}}

\newcommand{\Gm}{\mathbb{G}_m}

\renewcommand{\tilde}{\widetilde}
\newcommand{\Hom}{\operatorname{Hom}}

\newcommand{\Ext}{\operatorname{Ext}}

\newcommand{\Biext}{\operatorname{Biext}}

\newcommand{\Ker}{\operatorname{Ker}}

\renewcommand{\bar}[1]{\overline{#1}}

\def\spcheck{^\vee}

\newcommand{\et}{\mathrm{\acute{e}t}}

\newcommand{\cl}{Claim}
\theoremstyle{definition}

\newtheorem*{clm*}{\cl}

\theoremstyle{plain}

\def \fs {\mathrm {fs}}

\def \kfl {\mathrm {kfl}}

\def \overc#1{\overset {\lower 0.3ex \hbox{${\;}_{\circ}$}}{#1}}

\let\refsave=\ref
\def\ref#1{\textup{\refsave{#1}}}

\newcommand{\upc}{\overset{\circ}\to}
\newcommand\Cal{\mathcal}
\newcommand\define{\newcommand}
\renewcommand\bold{\Bbb}

\define\bZ{\bold Z}

\define\bQ{\bold Q}
\define\bN{\bold N}
\define{\cS}{\Cal S}
\define{\Lie}{\text{Lie}}
\define{\cH}{\Cal H}
\define{\cExt}{{\Cal E}xt}
\define{\cHom}{{\Cal H}om}
\define{\cO}{\Cal O}
\define{\an}{\mathrm{an}}

\def \upcf {\overset {\lower 0.3ex \hbox{${\;}_{\circ}$}} f}
\def \upcp {\overset {\lower 0.3ex \hbox{${\;}_{\circ}$}} p}
\def \upc#1{\overset {\lower 0.3ex \hbox{${\;}_{\circ}$}}{#1}}

\newcommand{\bs}{\backslash}
\newcommand{\Sig}{\Sigma}
\newcommand{\sig}{\sigma}

\newcommand{\cA}{\mathcal A}
\newcommand{\cB}{\mathcal B}
\newcommand{\cC}{\mathcal{C}}

\newcommand{\cP}{\mathcal{P}}
\newcommand{\cQ}{\mathcal Q}

\newcommand{\ep}{\varepsilon}

\newcommand{\ptpol}{\rm{ptpol}}
\newcommand{\pol}{\rm{pol}}

\begin{document}
\title
{Logarithmic abelian varieties, \\ Part VI: Local moduli and GAGF} 
\author{Takeshi Kajiwara, Kazuya Kato, and Chikara Nakayama}
\date{}
\maketitle
\centerline{\it{Dedicated to Professor Luc Illusie}} 
\setlength{\baselineskip}{1.0\baselineskip}
\begin{abstract}
\noindent 
This is Part VI of our series of papers on log abelian varieties.  In this part, we study local moduli and GAGF of log abelian varieties.
\end{abstract}
\section*{Contents}

\noindent \S\ref{sec:pol-tor}. GAGF for $\Gm$-, $\Gmlog$-, and $\Gmlog/\Gm$-torsors on weak log abelian varieties

\noindent \S\ref{s:ketpre}. K\'et presentation of a weak log abelian variety by a model, and equivalences with the categories of models

\noindent \S\ref{sec:fmod}. Moduli in the case of constant degeneration

\noindent \S\ref{s:dvr}. Weak log abelian varieties over complete discrete valuation rings, I

\noindent \S\ref{sec:GAGF1}. GAGF for log abelian varieties, I

\noindent \S\ref{sec:GAGF2}. GAGF for log abelian varieties, II

\noindent \S\ref{s:dvr2}. 
Weak log abelian varieties over complete discrete valuation rings, II

\section*{Introduction}
\renewcommand{\thefootnote}{\fnsymbol{footnote}}
\footnote[0]{Primary 14K10; 
Secondary 14J10, 14D06} 

  This is Part VI of our series of papers on log abelian varieties.
  In this part, we study the local moduli and the GAGF of log abelian varieties.
  In the next part, we will construct the global moduli of log abelian varieties, which is a main goal of this series of papers, by gluing local moduli by GAGF. 
  Thus the two main theorems of this paper are found in Sections \ref{sec:fmod} and \ref{sec:GAGF2}, that is, 
the description of local moduli (Theorem \ref{t:localmoduli}) and the GAGF for log abelian variety (Theorem \ref{thm5}). 

  In Part IV (\cite{KKN4}) of this series of papers, we needed several results on the category of weak log abelian varieties 
over a complete discrete valuation ring. 
  These were included in Part IV as a set of four sections titled 
\lq\lq Weak log abelian varieties over complete discrete valuation rings, I--IV.'' 
  In Sections \ref{s:dvr} and \ref{s:dvr2} of this paper, we prove some complementary results, 
which were announced in Part IV. 

  The results in Sections \ref{sec:pol-tor} and \ref{s:ketpre} themselves are important and used in the proof of the GAGF. 
  Section \ref{sec:GAGF1} proves the part of the full faithfulness of the GAGF.
  The theoretical dependence of the sections are as follows. 

$\ref{s:ketpre}\Rightarrow\ref{sec:GAGF1}$. 
\quad
$\ref{sec:pol-tor}, \ref{s:ketpre}, \ref{sec:fmod}, \ref{s:dvr}, \ref{sec:GAGF1} \Rightarrow \ref{sec:GAGF2}\Rightarrow \ref{s:dvr2}$. 

\smallskip

{\sc Acknowledgments.}
The authors thank Luc Illusie and Takeshi Saito for helpful discussions. 
The first author is partially supported by JSPS, Kakenhi (C) No.\ 24540035 and 
Kakenhi (C) No.\ 15K04811.
The second author is partially  supported by NFS grants DMS 1303421 and DMS 1601861.  
The third author is partially supported by JSPS, Kakenhi (C) No.\ 22540011, 
Kakenhi (B) No.\ 23340008, and Kakenhi (C) No.\ 16K05093.

\section{GAGF for $\Gm$-, $\Gmlog$-, and $\Gmlog/\Gm$-torsors on weak log abelian varieties}
\label{sec:pol-tor}
  In this section, we study GAGF for 
$\Gm$-torsors, $\Gmlog$-torsors, and $\Gmlog/\Gm$-torsors on a weak log abelian variety. 
  That is, we prove that these torsors on a weak log abelian variety over a complete 
strictly local noetherian ring are determined by the formal ones.

\begin{prop}\label{p:ff1}
  Let $S$ be an fs log scheme, and $X$ an fs log scheme over
  $S$. Assume that the underlying scheme $\overc{S}$ of $S$ is  the
  spectrum of a complete strictly local noetherian ring $(R,m)$, and that
  $\overc{X}$ is proper over $\overc{S}$. 
  For $n\geq 0$, let
  $S_n=\Spec(R/m^{n+1})$ (resp.\  $X_n=X\times_S S_n$) denote
  the fs log scheme with the inverse image log structure of $S$ (resp.\
  $X$).  
  Let $F=\Gm$ or $F=\Gmlog/\Gm$. 
Then the natural functors 
  $$ (\text{the groupoid of } F\text{-torsors on } X)
  \longrightarrow  (\text{the groupoid of } 
  F\text{-torsors on }\Hat{X})$$ 
  is fully faithful. Here an $F$-torsor on $\Hat{X}$ is 
  an inverse system $(G_n)_n$ of $F$-torsors $G_n$ on $X_n$
  such that $G_n$ is the pullback of $G_{n+1}$ by $X_n \to X_{n+1}$. 
\end{prop}
\begin{pf}
  For $F=\Gm$, it is enough to show the following (a) and (b): 

  (a)~$H^1(X,\Gm) \to 
\underset n \varprojlim\,  H^1(X_n, \Gm)$ is injective;

  (b)~$H^0(X,\Gm)\overset{\isom}{\longrightarrow} 
\underset n\varprojlim\,
  H^0(X_n, \Gm)$.

These are by Grothendieck's GAGF (\cite{EGA3}).
(The homomorphism in (a) is, in fact, an isomorphism.)

  Next, for $F=\Gmlog/\Gm$, the following (a$'$) and (b$'$) will suffice:

  (a$'$)~$H^1(X, \Gmlog/\Gm)\to H^1(X_n,\Gmlog/\Gm)$ is injective;

  (b$'$)~$H^0(X,\Gmlog/\Gm)
  \overset{\isom}{\longrightarrow}H^0(X_n,\Gmlog/\Gm)$. 

  To show them, we use the fact that $\Gmlog/\Gm$ on
  $(X_n)_{\et}$ is the inverse image of the one on $X_{\et}$ as
  abelian sheaves. 
(On the other hand, $\Gmlog$ (resp.\ $\Gm$) on $(X_n)_{\et}$
  is not necessarily 
the inverse image of the one on $X_{\et}$.) 
  By this fact, the above statements are deduced from the proper
  base change theorem in \'etale cohomology. 
\end{pf}

\begin{thm}\label{p:ffonLAV}
  Let the notation be as in the previous proposition. 
  Let $A$ be a weak log abelian variety over $S$ satisfying 
the conditions in $1.4.1$ of {\rm \cite{KKN5}}. 
  Let $F=\Gm, \Gmlog,$ or $\Gmlog/\Gm$. 
Then the natural functor
  $$ (\text{the groupoid of } F\text{-torsors on } A) \to 
  (\text{the groupoid of } F\text{-torsors on } \Hat{A})$$
  is fully faithful, equivalently,
  the homomorphism $H^i(A,F)\to \underset n\varprojlim\, H^i(A_n,F)$
  is injective (resp.\ an isomorphism) for $i=1$ (resp.\
  $i=0$). Here, $A_n$ is the pullback of $A$ to $S_n$, and 
an $F$-torsor on $\Hat{A}$ is 
  an inverse system $(G_n)_n$ of $F$-torsors $G_n$ on $A_n$
  such that $G_n$ is the pullback of $G_{n+1}$ by $A_n \to A_{n+1}$. 
\end{thm}

\begin{pf}
  Take a covering $P \to A$ by proper models as in Proposition 11.1 of \cite{KKN5}. 
  (Here we use the condition 1.4.1 of \cite{KKN5}.)  
  Note that, in this case, $P$ is the disjoint union of proper models over $S$ because 
$R$ is strictly local.
  Let $\cC(-)$ denote the groupoid of $F$-torsors. 
  Consider the following essentially commutative diagram:
$$\begin{matrix} \cC(A)&  \to& \cC(P) &\rightrightarrows  & \cC(P\times_A P)\\
\downarrow && \downarrow&& \downarrow \\
\underset n\varprojlim\, \cC(A\times_S S_n)&  \to& \underset n\varprojlim\, \cC(P\times_S S_n) &\rightrightarrows  
&\underset n\varprojlim\, \cC((P\times_A P)\times_S S_n).
\end{matrix}$$

  Since we have a descent for each horizontal row, the full faithfulness for $A$ 
is reduced to those for $P$ and for $P \times_AP$. 
  Hence, if $F=\Gm$ or $\Gmlog/\Gm$, the desired full faithfulness is reduced to 
the previous proposition. 

  Next consider the following commutative diagram: 
  \begin{multline*}
  \begin{CD}
    0 @>>> H^0(\Gm)@>>> H^0(\Gmlog) @>>> H^0(\Gmlog/\Gm) 
    @>>>  \\
    @. @VV{\vert\wr}V @VVV @VV{\vert\wr}V @. \\
    0 @>>> \underset{n}{\varprojlim}\, H^0(\Gm)_n @>>>
    \underset{n}{\varprojlim}\, H^0(\Gmlog)_n @>>>
    \underset{n}{\varprojlim}\, H^0(\Gmlog/\Gm)_n @>>>
  \end{CD} \\
  \begin{CD}
    H^1(\Gm) @>>> H^1(\Gmlog) @>>> H^1(\Gmlog/\Gm) \\
    @V{\cap}VV @VVV @V{\cap}VV \\
    \underset{n}{\varprojlim}\, H^1(\Gm)_n @>>>
    \underset{n}{\varprojlim}\, H^1(\Gmlog)_n @>>>
    \underset{n}{\varprojlim}\,  H^1(\Gmlog/\Gm)_n.
  \end{CD}
  \end{multline*}
  Here $H^*(\Gm)$ (resp.\ $H^*(\Gm)_n$) etc.\ denote $H^*(A,\Gm)$
  (resp.\ $H^*(A_n, \Gm)$) etc..
Note that the upper row is exact, and
  the lower one is a complex. 
  To reduce the full faithfullness in the case where $F=\Gmlog$ to 
the cases where $F=\Gm$ and $F=\Gmlog/\Gm$, 
it is enough to see that the lower row is
  exact at $\underset n \varprojlim\, H^0(\Gm)_n$, 
$\underset n\varprojlim\, H^0(\Gmlog)_n$,
  and $\underset n\varprojlim\, H^1(\Gm)_n$. 
  Among these, the only nontrivial one is the 
exactness at $\underset n\varprojlim\, H^1(\Gm)_n$.
  This is proved by the fact that the inverse system of 
$H^0(A_n,\Gmlog)=H^0(S_n,\Gmlog)$ (Proposition 6.2 of \cite{KKN5}) satisfies the Mittag-Leffler condition. 
\end{pf}

\begin{prop}
  Let $S$, $A$ be as in Proposition \ref{p:ffonLAV}. 
  Let $F=\Gm, \Gmlog,$ or $\Gmlog/\Gm$. 
Then 
we have a Cartesian diagram
  $$
  \begin{CD}
    \Ext(A,F) @>>> H^1(A,F) \\
    @VVV @VVV \\
    \underset{n}{\varprojlim} \Ext(A_n,F) @>>> 
    \underset{n}{\varprojlim} H^1(A_n,F),
  \end{CD}$$
where the all arrows are injective. 
\end{prop}
\begin{pf}
  This follows from Proposition \ref{p:ffonLAV} and \cite{KKN5} 
  Lemma 12.3.
\end{pf}

\section{K\'et presentation of a weak log abelian variety by a model, and equivalences with the categories of models }\label{s:ketpre}
  To prove GAGF for log abelian varieties, it is crucial to replace log abelian varieties by models, for log abelian varieties are only functors but models are fs log schemes, and we can apply classical theories to models. 
  To this end, we establish the category equivalence between the category $\cA$ of weak log abelian varieties (with additional data) and the category $\cB$ of models (with additional data) in this section.
  There are many variants of this ``recovering from models''-type statements. 
  We already mention one of them in \cite{KKN5} Remark 11.9.
  Another statement is \cite{KKN3} Proposition 4.4, which is a part of 1-dimensional case of Theorem \ref{t:A,B} below but contains a mistake.  
  See Remark \ref{r:la3mistake} (1)  below.  

\begin{para} 
\label{setting}
  Since the equivalence discussed here is based on the k\'et presentation of a weak log abelian variety by a model, we first describe it. 
  In Section 11 of \cite{KKN5}, we covered a weak log abelian variety by models with respect to \'etale topology.
  The point here is that if we work with k\'et topology instead of \'etale topology, we can make a much simpler cover, and even the relation can be simply described.

  Let $A$ be a weak log abelian variety over an fs log scheme $S$. 
  Let $G$ be its semiablian part (\cite{KKN2} 4.4, \cite{KKN4} 1.7). 
  Assume that there are an admissible pairing (\cite{KKN2} 7.1)
$$X \times Y \to \Gmlog/\Gm$$ 
on $S$, where $X$ and $Y$ are finitely generated free $\bZ$-modules, and an isomorphism 
$$A/G \cong \cHom(X, \Gmlog/\Gm)^{(Y)}/\overline Y$$
(see \cite{KKN4} 1.3 for the definition of the right-hand-side). 
  Let $\tilde A$ be the fiber product of $A \to \cHom(X, \Gmlog/\Gm)^{(Y)}/\overline Y \leftarrow \cHom(X, \Gmlog/\Gm)^{(Y)}$.

  Assume further that there are a homomorphism $\cS \to M_S/\cO^{\times}_S$ from an fs monoid and an $\cS$-admissible pairing 
$$X \times Y \to \cS^{\gp}$$ 
which lifts the above $\Gmlog/\Gm$-valued pairing. 
  Let $$C=\{(N,l)\in \Hom(\cS, \bN) \times \Hom(X,\bZ)\,|\,l(X_{\Ker(N)})=0\}$$
(cf.\ \cite{KKN1} 3.4.2). 
\end{para}

\begin{para}
  Let $\Sig$ be a $Y$-stable fan in $C$ (\cite{KKN4} 2.6). 

  We review the definition that $\Sig$ is wide (\cite{KKN5} 10.1). 
  Let $\sig$ be a {\it cone}, that is, a finitely generated $\bQ_{\ge0}$-submonoid of $C_{\bQ_{\ge0}}$. 
  We say that $\sig$ is {\it wide} if for any $(N,l) \in C$, we have $(N, \ep l)\in \sig$ for any $\ep \in \bQ$ such that $|\ep|$ is sufficiently small. 
  We say that $\Sig$ is {\it wide} if it owes a wide cone.
  Note that if $\Sig$ is wide, a wide cone in $\Sig$ is unique. 

  The result in this section roughly says that if $\Sig$ is complete and wide, $A$ recovers from its $\Sig$-part (or $\Sig$-model) $A^{(\Sig)}$. 
  (See \cite{KKN4} 3.1 for the definition of completeness.)

  Note that the first standard fan and the second standard fan (cf.\ \cite{KKN5} 4.6) are complete and wide.
  Their wideness is by \cite{KKN5} Proposition 10.3.
\end{para}

\begin{para}
\label{ketpre}
  We have the k\'et presentation of $\tilde A$ and then of $A$ as follows. 
  Assume that there is a prime number $\ell$ which is invertible on $S$. 
  We fix such an $\ell$. 

  Let $$Z:=\coprod_{n\geq 0} {\tilde A}^{(\Sig)} = {\tilde A}^{(\Sig)} \times \bN$$ 
and consider the morphism 
$$Z \to \tilde A\;;\;(x,n)\mapsto x^{\ell^n}.$$ 
It is k\'et surjective, that is,  surjective with respect to the k\'et topology because the induced 
$\coprod_{n\geq 0}{\tilde A}^{(\sig)} \to A$ is already k\'et surjective for a wide cone $\sig$. 

  Let $R:=Z\times_A Z$. 
  It is the disjoint union of the $(m,n)$-part. Assume $m\geq n$. Then the $(m,n)$-part of $R$ is
isomorphic to $\tilde A^{(\Sig(m,n))}\times G[\ell^n]$, where $\Sig(m,n)=\Sig \sqcap \ell^{n-m}\Sig$ ($\ell^{n-m}\Sig:=\{\ell^{n-m}\sig\;|\;\sig \in \Sig\}$), and the map to the fiber product is 
$(a,b)\mapsto (a, a^{\ell^{m-n}}b)$. 
  Here we recall that for fans $\Sig_1$ and $\Sig_2$, the fan $\Sig_1 \sqcap \Sig_2$ is defined as 
$\{\sig_1 \cap \sig_2 \;|\; \sig_1 \in \Sig_1, \sig_2 \in \Sig_2\}$ (\cite{KKN1} Definition 5.2.15). 
  If $m \leq n$, the $(m,n)$-part of $R$ is 
isomorphic to $\tilde A^{(\Sig(n,m))}\times G[\ell^m]$, and the map to the fiber product is 
$(a,b)\mapsto (a^{\ell^{n-m}}b,a)$.

  Next the group structure of $\tilde A$ is characterized by the partial group law of $Z$ as follows. 
  Below, in general, a {\it partial group law} on a sheaf $F$ simply means a map from a subsheaf of $F \times F$ to $F$ (we do not impose associativity etc.). 
  Let $(\tilde A^{(\Sig)}\times \tilde A^{(\Sig)})'$ be the subsheaf of $\tilde A^{(\Sig)}\times \tilde A^{(\Sig)}$ consisting of the sections $(x,y)$ such that 
the product of the images of $x$ and $y$ in $\cHom(X,\Gmlog/\Gm)$ belongs to the $\Sig$-part. 
  Then we have a partial group law 
$(\tilde A^{(\Sig)}\times \tilde A^{(\Sig)})' \to \tilde A^{(\Sig)}$
on $\tilde A^{(\Sig)}$ by restricting the group law of $\tilde A$. 
  From this, we give $Z$ a partial group law 
defined by 
$$(x, n)(y, n) = (xy, n)  \qquad (x, y) \in (\tilde A^{(\Sig)}\times \tilde A^{(\Sig)})', n \in \bN.$$ 
  Then the group structure of $\tilde A$ is characterized by the unique one which is compatible with this partial group law of $Z$.

  Finally, we recover $A$ by dividing by $\bar Y$ as follows. 
  Since there exists a wide cone in $\Sig$, $\tilde A^{(\Sig)}$ contains the 0-section of $\overline Y$.  
  Since $\Sig$ is $Y$-stable, $\tilde A^{(\Sig)}$ contains the other sections of $\overline Y$.
  Thus there is a map $\overline Y \to Z$ via the 0-th component and we recover the homomorphism $\overline Y \to \tilde A$, and we recover $A=\tilde A/\bar Y$.
\end{para}

  Taking the above observation into account, we introduce the following categories $\cA$ and $\cB$. 

\begin{para}\label{A,B}
  We define the categories $\cA$ and $\cB$.
  First we define the objects of these categories. 

  Let $S$ be an fs log scheme, and $\ell$ a prime number which is invertible 
on $S$. 

  The objects of $\cA$ are pairs $(A, \Sig)$, where $A$ is a weak log abelian variety over $S$ and $\Sig$ is a subsheaf 
of $Q:=\cH om(\bar X, \Gmlog/\Gm)^{(\bar Y)}/\bar Y$ 
coming from complete and wide fans.
  Here $\cH om(\bar X, \Gmlog/\Gm)^{(\bar Y)}$ is the associated one to the admissible pairing $\bar X \times \bar Y \to \Gmlog/\Gm$ determined by $A$. 
  That a subsheaf $\Sig$ comes from complete and wide fans means that \'etale locally on $S$, there are the data $X \times Y \to \cS^{\gp}$ and 
$\cS \to M_S/\cO_S^{\times}$ as in \ref{setting} 
inducing the admissible pairing $\bar X \times \bar Y \to \Gmlog/\Gm$ determined by $A$, and a complete and wide fan $\Sig'$ in the associated $C$ 
such that $\Sig$ coincides with 
$\cH om(\bar X, \Gmlog/\Gm)^{(\Sig')}/\bar Y$. 

  The objects of $\cB$ are 5-ples $(P, e, G, Q, \Sig)$, where 
$G$ is a
semiabelian scheme over $S$, 
$Q$ is a sheaf of abelian groups coming from admissible pairings, 
$\Sig$ is a subsheaf of $Q$ coming from complete and wide fans, 
$P$ is a $G$-torsor over $\Sig$ 
endowed with a partial group law $(P\times P)' \to P$, and 
represented by a log algebraic space in the first sense over $S$ (\cite{KKN4} 10.1),  
and $e$ is a section of $P$ over $S$ which maps to the unity of $Q$
($e$ is called the origin of $P$), 
satisfying the conditions 1, 2 and 3 below.
  Here that $Q$ comes from admissible pairings means that, \'etale locally on $S$, 
there are the data $X \times Y \to \cS^{\gp}$ and 
$\cS \to M_S/\cO_S^{\times}$ as in \ref{setting}
such that $Q$ is isomorphic to 
$\cH om(\bar X, \Gmlog/\Gm)^{(\bar Y)}/\bar Y$.
  That $\Sig$ comes from complete and wide fans means that \'etale locally on $S$, there are the same data as in \ref{setting} and a complete and wide fan $\Sig'$ in the associated $C$ such that there is an isomorphism between $Q$ and 
$\cH om(\bar X, \Gmlog/\Gm)^{(\bar Y)}/\bar Y$ via which 
$\Sig$ coincides with $\cH om(\bar X, \Gmlog/\Gm)^{(\Sig')}/\bar Y$. 
  The $(P\times P)'$ is the inverse image of $\Sig$ by $P \times P \to Q \times Q \to Q; (x,y) \mapsto (x,y) \mapsto xy$.

The three conditions for objects of $\cB$ are as follows:  

1. The partial group law on $P$ is compatible with the action of $G$ on $P$ and with $(x,y)\mapsto xy$ on the quotient $Q$.

2. The partial group law on $P$ satisfies the following three conditions (a)--(c).  
  We say that $xy$ is {\it defined} if $(x,y) \in P \times P$ belongs to $(P \times P)'$.

\quad $(\rm a)$ For any $x,y,z \in P$, if $xy, (xy)z, yz$ are defined, then $x(yz)$ is also defined and $(xy)z=x(yz)$. 

\quad $(\rm b)$ For any $x \in P$, $ex$ is defined and $ex=x$. 

\quad $(\rm c)$ For any $x,y\in P$, if $xy$ is defined, then $yx$ is also defined and $xy=yx$. 

3. $P$ is separated over $S$.
\end{para}

\begin{rem}
\label{r:la3mistake}
(1)
  The category $\cB$ defined in \cite{KKN3} Proposition 4.4 is essentially the same as the full subcategory consisting of 
1-dimensional objects of the category $\cB$ here. 
  But the condition 

(h) The sheaf $W$ is quasi-separated over $S$. 

\noindent in \cite{KKN3} Proposition 4.4, which corresponds to the condition 3 in the above, should be corrected as follows.

(h) The sheaf $W$ is separated over $S$. 

\noindent Without this change, Proposition 4.4 of \cite{KKN3} is not valid.
  In fact, in the last part of the proof, $E$ is proved to be quasi-separated but it does not imply that $E$ is separated (see \cite{KKN2} 11.6 for a counter example).  
  After changing the condition (h) as above, we can prove that $E$ is separated by the argument with the 0-section in the last part of the proof of Theorem \ref{t:A,B} below so that Proposition 4.4 of \cite{KKN3} becomes valid and it is not necessary to change the remaining part of \cite{KKN3}. 

(2) There is a similar mistake in \cite{KKN4}.
  See Remark \ref{r:la4mistake} below. 

(3) In the condition 2 (a) and (c), that $x(yz)$ and $yx$ are defined is a conclusion of the latter half of the condition 1. 
  But in \cite{KKN3} 4.3, the condition (e), which corresponds to the latter half of the condition 1, was introduced too late and 
it is not automatic that $x(yz)$ and $yx$ are defined.
  Hence, the conditions (a), (c) in \cite{KKN3} 4.3 should be replaced by the conditions (a), (c) in 2.4 in this paper. 
\end{rem}

\begin{para} (Continuation of the definitions of $\cA$ and $\cB$.)
  We define the morphisms of the categories $\cA$ and $\cB$.

Morphisms in $\cA$. $(A, \Sig)\to (A', \Sig')$ is a homomorphism $A\to A'$ which induces $\Sig \to \Sig'$.

Morphisms in $\cB$. $(P, e, G, Q, \Sig) \to (P', e', G', Q', \Sig')$ is a triple $(a,b,c)$ consisting of a homomorphism $a\colon G\to G'$, a homomorphism 
$b\colon Q\to Q'$ which sends $\Sig$ into $\Sig'$, and a morphism $c\colon 
P\to P'$ which is compatible with $a$ and with $b$ and which 
commutes with the partial group laws. 
\end{para}

\begin{thm}\label{t:A,B}
  The natural functor 
$$\cA\overset{\simeq}\to \cB\;;\; (A,\Sig)\mapsto (A^{(\Sig)},e,G,A/G,\Sig)$$ 
gives an equivalence of categories. 
  Here $G$ is the semiabelian part of $A$. 
\end{thm}

\begin{pf}
  We give the inverse functor $\cB \to \cA$. 
  Let $(P, e, G, Q, \Sig)$ be an object of $\cB$.
  To recover $A$, we use the k\'et presentation of $A$ in \ref{ketpre}. 
  
  First note that by Theorem 7.6 in \cite{KKN2}, $Q$ decides the sheaves $\overline X$, $\overline Y$ and the pairing $\overline X \times \overline Y \to \Gmlog/\Gm$ such that 
$Q\cong \cHom(\overline X, \Gmlog/\Gm)^{(\overline Y)}/\overline Y$ globally.
  Define 
$$\tilde P 
=P \times_Q\tilde Q,$$
where 
$$\tilde Q:=\cHom(\overline X, \Gmlog/\Gm)^{(\overline Y)}.$$
  The partial group law on $P$ and the group law of 
$\tilde Q$ induce a partial group law on $\tilde P$.

  By \ref{ketpre}, we can recover $\tilde A$ with the group structure as follows. 

  Let $Z:=\coprod_{n\geq 0} {\tilde P}= {\tilde P} \times \bN$. 

  Let $R:=\coprod_{m\geq n\geq 0} \tilde P^{(\Sig(m,n))}\times G[\ell^n] \sqcup
\coprod_{n> m\geq 0} \tilde P^{(\Sig(n,m))}\times G[\ell^m]$. 

  We define the morphism 
$R \to Z\times_SZ$ as the induced one by the morphisms
$$\tilde P^{(\Sig(m,n))}\times G[\ell^n] \to Z\times_SZ; 
(a,b)\mapsto ((a,m), (a^{\ell^{m-n}}b,n)) \quad (m \geq n)$$ and 
$$\tilde P^{(\Sig(n,m))}\times G[\ell^m] \to Z\times_SZ; 
(a,b)\mapsto ((a^{\ell^{n-m}}b,m),(a,n)) \quad (n > m),$$ where $\tilde P^{(\Sig(m,n))}$ is the pullback of the $\Sig(m,n)$-part of $Q$ by $\tilde P \to Q$. 

  Let $\tilde A$ be the k\'et difference cokernel of $R \rightrightarrows Z$. 

  Let $Z \to \tilde Q$ be the morphism defined on the $n$-th component of $Z$ by $\tilde P \to \tilde Q \overset {\ell^n} \to \tilde Q$, 
which is surjective since $\Sig$ is wide.
  This morphism induces a surjection $\tilde A \to \tilde Q$.

  Define the group structure of $\tilde A$ as follows. 
  Let $(x,m)$, $(y,n)$ ($x,y \in \tilde P, m,n \geq 0$) be in $Z$. 
  Since $\Sig$ is wide, k\'et locally, there is a sufficiently big $n'\geq 0$ such that 
$\overline{(x,m)}=\overline{(x',n')}$ and $\overline{(y,n)}=\overline{(y',n')}$ for some $(x', y') 
\in (\tilde P \times \tilde P)'$, where $\overline a$ is the image of $a \in Z$ in $\tilde A$.
  Then we define $\overline{(x,m)}\cdot \overline{(y,n)}=\overline{(x'y', n')}$, which does not depend on the choices and gives a group law on $\tilde A$ by the condition 2. 
  The morphism $\tilde A \to \tilde Q$ is a homomorphism.
  The morphism $\tilde P \to Z \to \tilde A$, where the first morphism is the $0$-th inclusion, is injective, and we identify $\tilde P$ with the image of this injection.  
  Then $\tilde P$ coincides with the inverse image of $\Sig$ in $\tilde A$ by the homomorphism $\tilde A \to \tilde Q \to Q$.

  Define $\overline Y\to \tilde A$ as follows. 
  Let $y \in \overline Y$. 
  K\'et locally, there is an $n \ge0$ and a section $y_1$ of the $\Sig$-part of 
$\tilde Q$ such that $y_1^{\ell^n}$ coincides with the image of $y$ in $\tilde Q$. 
  We may assume that the image of $y_1$ in $Q$ which belongs to $\Sig$ comes from a section $a$ of $P$. 
  Then $a$ defines a lift $\tilde y_1$ to $\tilde P$ of $y_1$.
  Then we define $f(y)=\overline{(\tilde y_1, n)} \in \tilde A$. 
  This is independent of choices and defines an injective homomorphism $f\colon \overline Y\to \tilde A$.

  Define $A=\tilde A/\overline Y$. 
  Then we have an injective morphism $P \to A$ 
and a homomorphism $A \to Q$, and we can identify $P$ with the $\Sig$-part of $A$. 

  We prove that $(A, \Sig)$ is an object of $\cA$. 

  First, we consider on $G$. 
  The composite of the inclusion $G \to P$ via $e$ and the inclusion $P \to A$ is an injective homomorphism. 
  Then a direct calculation shows that $A/G$ is naturally isomorphic to $Q$. 
  (By the arguments in \cite{KKN2}, 9.2 and 9.3, $G$, $\overline X$, and $\overline Y$ are those determined by $A$.)
  Hence, the second condition 1.6 (2) of \cite{KKN4} in the definition of weak log abelian variety is satisfied. 

  A big nontrivial point is that any fiber of $A$ comes from an admissible and nondegenerate log 1-motif, that is, that the first condition (1) in Definition 1.6 of \cite{KKN4} in the definition of weak log abelian variety is satisfied. 
  The proof of this point is as follows. 
  We may and do consider the case of constant degeneration. 
  We assume that $\overline X = X$ and $\overline Y = Y$. 
  Let $T=\cHom(X,\Gm)$.
  Then we have $\tilde Q=(T_{\log}/T)^{(Y)}$. 

  Using the last statement of \cite{KKN2} Theorem 7.3 (1), we see that 
the exact sequence
$0\to G\to \tilde A \to (T_{\log}/T)^{(Y)}\to 0$ is obtained by the pushout of the exact sequence 
$0\to T\to T_{\log}^{(Y)} \to (T_{\log}/T)^{(Y)}\to 0$ by some homomorphism $h\colon T\to G$.
  Then $\tilde P=\tilde A^{(\Sig)}$ is the pushout of 
$T_{\log}^{(\Sig)} \leftarrow T \overset h \to G$. 
  We prove that the representability of $\tilde P$ implies that 
$h$ is an isomorphism to the torus part of $G$ so that $A$ comes from an admissible and nondegenerate log 1-motif.

  Let $H$ be the kernel of $h$. 
  If $h$ is not an isogeny to the torus part of $G$, the nonrepresentable $H_{\log}^{(\Sig)}/H$ survives in the pushout $\tilde P$, 
where $H_{\log}^{(\Sig)}=H_{\log}\cap T_{\log}^{(\Sig)}$.
  Hence the pushout is not representable. 
  More precisely, there is an injective morphism $H_{\log}^{(\Sig)}/H \to \tilde P$. 
  On the other hand, there is a wide cone $\sig$ in $\Sig$, and, by Lemma 9.10 of \cite{KKN4}, the induced morphism 
$(H_{\log}\cap T_{\log}^{(\sig)})/H \to \tilde P$ should be constant, which is a contradiction. 
  Hence, $h$ is an isogeny to the torus part of $G$.

  Next, if $h$ is not an isomorphism to the torus part of $G$, then 
the quotient $P=Y\bs  G_{\log}^{(\Sigma)}/H$ as a sheaf for the \'etale topology is not representable, which is a contradiction.  
  We prove this.  
  First, the kfl (kummer log flat) quotient $P^{\mathrm{kfl}}=
(Y\bs  G_{\log}^{(\Sigma)}/H)^{\mathrm{kfl}}$ is representable because 
$Y\bs  G_{\log}^{(Y)}$ is a weak log abelian vareity, 
the kfl quotient $((Y\bs G_{\log}^{(Y)})/H)^{\mathrm{kfl}}$ is also a weak log abelian variety, and $P^{\mathrm{kfl}}$ is a model of this.
  But there is a difference between these two kinds of quotients $P$ and $P^{\mathrm{kfl}}$.
  To see it, let $X' \subset X$ be a nontrivial subgroup of finite index.  
  Then the map $\cHom(X,\Gmlog) \to \cHom(X',\Gmlog)$ is kfl surjective but not \'etale surjective.  This is reduced to the case $X=\bZ$.
  Thus $P$ is not a sheaf for the kfl topology.
  By \cite{Kato:FI2} Theorem 3.1, $P$ is not representable, a contradiction. 

  The rest is the separability, that is, the third condition 1.6 (3) of \cite{KKN4} in the definition of weak log abelian variety.
  By the next Lemma \ref{l:sep}, it is enough to show that the $0$-section $0\colon S \to A$ is represented by proper morphisms. 
  But this $0$ factors as $S \overset e \to P \overset i \to A$.
  By the condition 3, $e$ is represented by finite morphisms. 
  To prove that $i$ is represented by proper morphisms, we cover $A$ by $\tilde A^{(\sig)}$ with $\sig$ being a various cone. 
  Then the base-changed morphism of $i$ is $\tilde A^{(\Sig\sqcap \sig)} \to \tilde A^{(\sig)}$, which is a log blow-up. 
  Hence $i$ is represented by proper morphisms. 
  Therefore, $0$ is represented by proper morphisms. 
\end{pf}

\begin{lem}
\label{l:sep}
  In the definition {\rm \cite{KKN2}} Definition $4.1$ of log abelian variety and in the definition {\rm \cite{KKN4}} Definition $1.6$ of weak log abelian variety, the third condition 

$(*)$ The diagonal morphism $A \to A \times_S A$ is represented by finite morphisms.

\noindent 
can be replaced by the following condition{\rm:} 

$(**)$ The $0$-section $S \to A$ is represented by proper morphisms.
\end{lem}

\begin{pf}
  It is enough to show that the condition $(**)$ implies that the $0$-section is represented by finite morphisms (cf.\ the last paragraph of \cite{KKN2} 4.2). 
  Hence it suffices to see that the $0$-section is always represented by morphisms with finite fibers. 
  This is reduced to the case with constant degeneration. 
  We use the notation as in \ref{setting}. 
  Cover $A$ with $\tilde A^{(\sig)}$ with $\sig$ being a various cone. 
  Then the base-changed morphism of the $0$-section $\overline Y \cap \tilde A^{(\sig)} \to \tilde A^{(\sig)}$ is a morphism with finite fibers. 
  Hence we conclude that the $0$-section is represented by morphisms with finite fibers. 
\end{pf}

\begin{rem}
\label{r:la4mistake}
  We correct the related part of \cite{KKN4}.   
  In the last paragraph of \cite{KKN4} 9.12, which is a part of the proof of Proposition 9.2 of \cite{KKN4}, is wrong.
  More precisely, in the notation there, the finiteness of $I \cap \overline Y \to I$ does not imply that of $\overline Y \to \tilde L$ (cf.\ Remark \ref{r:la3mistake} (1)). 

  We modify the argument as follows. 
  Let the notation be as in there. 
  We may assume that each $S_{\lambda}$ is of finite type over $\bZ$. 
  Take a complete fan $\Sigma$ by \cite{KKN1} Theorem 5.2.1. 
  Then the base-changed morphism 
$$0'\colon S^{(\Sig)} \to L^{(\Sig)}$$ 
of the $0$-section $S \to L$ is represented by finite morphisms, where $S^{(\Sig)}=S \times_LL^{(\Sig)}$. 
  Since each $S_{\lambda}$ is of finite type over $\bZ$, by the part of Theorem 8.1 in \cite{KKN4} proved till there, $L^{(\Sig)}$ over $S_{\lambda}$ is represented by an algebraic space with fs log structure of finite presentation. 
  Note that $S^{(\Sig)}$ is a log blow-up of $S$. 
  Hence $0'$ comes from a finite morphism over some $S_{\lambda}$.
  On the other hand, by the argument in the last paragraph of the proof of Theorem \ref{t:A,B}, 
the morphism $i\colon L^{(\Sig)} \to L$ is represented by proper morphisms and 
we see that the composite $i\circ 0'\colon S^{(\Sig)} \to L$ over $S_{\lambda}$ is represented by proper morphisms.

  Since $S^{(\Sig)} \to S$ is proper, if the $0$-section $S_{\lambda} \to L$ over $S_{\lambda}$ is represented by morphisms locally of finite type, it is also represented by proper morphisms.
  Cover $L$ by $\tilde L^{(\sig)}$ with a various cone $\sigma$, and the base-changed $0$-section is $\tilde L^{(\sig)} \cap \overline Y \to \tilde L^{(\sig)}$.
  Since $\tilde L^{(\sig)} \cap \overline Y \to S_{\lambda}$ is represented by a morphism locally of finite type, the base-changed $0$-section is also. 
  Thus $S_{\lambda} \to L$ is represented by proper morphisms. 
  By Lemma \ref{l:sep}, we conclude that $L$ satisfies the separability over $S_{\lambda}$. 
  This completes the proof of Proposition 9.2 of \cite{KKN4}. 
\end{rem}

By Theorem \ref{t:A,B}, we have a new definition of a weak log abelian variety with a wide fan. 

\begin{cor}
  Let $S$ be an fs log scheme. 
  An abelian sheaf $A$ on $(\fs/S)_{\et}$ is a weak log abelian variety if and only if 
there are an admissible pairing $\bar X \times \bar Y \to \Gmlog/\Gm$,  a semiabelian scheme $G$, an exact sequence $0\to G\to A \to \cHom(\bar X, \Gmlog/\Gm)^{(\bar Y)}/\bar Y\to 0$, and the following is satisfied. 
  \'Etale locally on the base, there is a complete and wide fan $\Sig$ such that $A^{(\Sigma)}$ is represented by a proper log algebraic space in the first sense.
\end{cor}

\begin{pf} 
The data give an object of $\cB$. 
\end{pf}

\section{Moduli in the case of constant degeneration}
\label{sec:fmod}
  In this section, we study the local moduli space of principally polarized log abelian varieties with constant degeneration and level structure.
  Let us fix a free abelian group $Y$ of rank $r$ throughout this section. Let $S_{\Q}(Y)$ be the set of symmetric bilinear forms $Y \times Y \to \Q$. Let $\sig \subset S_{\Q}(Y)$ be a finitely generated $\Q$-subcone satisfying the following condition:
$b(y, y) \geq 0$ for any $b \in \sig$ and $y \in Y$. 

\begin{defn} 
  Let $S$ be an fs log scheme. 
  A polarized log 1-motif (\cite{KKN2} Definition 2.2, \cite{KKN2} Definition 2.8) $M = [Y \to G_{\log}]$ over $S$ of type $(X,Y)$ is {\it of degeneration along $\sig$} if, for any $s \in  S$ and any homomorphism $\nu\colon 
(M_S/{\cal O}_S^\times)_{\bar s}\to \N$, the
composite map
$$Y\times Y =Y_{\bar s}\times Y_{\bar s} \overset{\phi\times \mathrm{id}}\to X_{\bar s}\times Y_{\bar s} \overset{\langle\;,\;\rangle}\to (M_S/\cO^\times_S)^{\gp}_{\bar s} \overset{\nu^{\gp}}\to \Z$$
is in the cone $\sig$. Here 
$\phi: Y \to X$ is the homomorphism induced by the polarization. 
  We say that a polarized log abelian variety with constant degeneration of type $(X, Y)$ is {\it of degeneration along $\sig$} if the corresponding polarized log $1$-motif is so.
\end{defn}

\begin{para}
  Let $g \ge r$ and $n\ge 3$ be integers, and $S$ an fs log scheme over $\Z[1/n]$. We define a functor $F = F_{g,r,n,\sig} : (\fs/S) \to$ (set) as follows: for each object $U$ of $(\fs/S)$,
$$F(U) :=$$ 

\noindent $\{$principally polarized log 1-motif of degeneration along $\sig$ over $U$
of the form  $[Y\to  (-)_{\log}]$ whose abelian part is of dimension $g - r$ 
and endowed 
with $n$-level structure$\}/\cong$. 

\end{para}

We remark that, by Theorem 3.4 in \cite{KKN2}, the functor $F$ is the moduli functor of principally polarized log abelian varieties of type $(*, Y )$ and of degeneration along $\sig$ whose abelian quotient of the semiabelian part is of dimension $g-r$ 
and endowed 
with $n$-level structure.

\begin{para} 
  We prove that the above functor $F = F_{g,r,n,\sig}$ is pro-represented by an fs log formal scheme.

The fs log formal scheme which represents $F$ is described as follows. 
  Let $Z$ be the moduli space of principally polarized abelian varieties of dimension $g - r$ with $n$-level structure over $\overset {\circ} S$-scheme, where $\overset {\circ} S$ is the underlying scheme of $S$. 
  The reason why we impose $n \geq  3$ is that the fine moduli space of the above moduli space exists. 
  Let $B$ be a universal abelian scheme over $Z$, and $\tilde Z := \cHom(Y,B^*)$. Hence, $\tilde Z$ is, \'etale locally on $S^{\circ}$, isomorphic to the $r$-times fiber product of $B^*$ over $Z$. This $\tilde Z$ represents a moduli functor of extensions of $B$ by the torus $T = \cHom(Y,\Gm)$ because $$\cExt(B,T) = \cHom(Y,\cExt(B,\Gm)) = \cHom(Y,B^*).$$
\end{para}

\begin{thm}
\label{t:localmoduli}
The functor $F$ is pro-represented by an fs log formal scheme which is, \'etale locally on $\tilde Z$, isomorphic to $(\Spec(\cO_{\tilde Z}[\sig{\spcheck}]/(I^n)))_n$. 
Here 
$$\sig{\spcheck} := \{x \in \mathrm{Sym}^2_{\Z}(Y )\;|\;  b(x) \geq  0 \text{ for every }b \in \sig\},$$  
$I$ is the ideal of the semigroup ring $\cO_{\tilde Z}[\sig{\spcheck}]$ 
generated by $y\otimes y\in \sig{\spcheck}$ for all $y\in Y -\{0\}$ and this formal scheme is endowed with the fs log structure associated to $M_S \times \sig{\spcheck}$.
\end{thm}

\begin{pf}
  We work over $\tilde Z$. 
  Let $0 \to T \to G \to B \to 0$ be the universal extension over $\tilde Z$.
  We have to parametrize a subset of the set of the homomorphism $Y \to G_{\log}$ such that 
$Y \to G_{\log} \to B$ is the given universal $Y \to B^* \cong B$. 
  We may assume that $F$ is not empty. 
  Then there is a bijection between the above set and the set $\Hom(Y,T_{\log})=\Hom(Y\otimes Y,\Gmlog)$.
  The condition that $Y \to G_{\log}$ is a principally polarized log $1$-motif of degeneration along $\sig$ is equivalent to the condition that the corresponding $Y \otimes Y \to \Gmlog$ comes from a section to 
$\Spec(\cO_{\tilde Z}[\sig{\spcheck}]/(I^n))$ for some $n$. 
\end{pf}

\begin{rem} The following three conditions are equivalent. (i) $F$ is not empty. (ii) If $b$ is in the interior of $\sig$ and $ y\in Y -\{0\}$, then $b(y,y)>0$. (iii) The subset $\{y \otimes y\;|\; y\in Y -\{0\} \}$ of $\sig{\spcheck}$ does not contain any invertible element of $\sig{\spcheck}$.
\end{rem}

\section{Weak log abelian varieties over complete discrete valuation rings, I}
\label{s:dvr}
  In a former part \cite{KKN4} of this series of papers, we study the category of weak log abelian varieties 
over a complete discrete valuation ring. 
  Though we proved there only what are necessary for the existence of proper models, we 
announced some related results mainly concerning the polarizations, which we prove here. 

\begin{para}
  First we recall the situation. 
  Let $K$ be a complete discrete valuation field with valuation ring $O_K$. 

  Fix an fs log structure $N$ on $\eta:=\Spec(K)$ charted by the stalk $\cS$ of $N/\cO^\times_{\eta}$,
and let $M$ be its direct image log structure on $\Spec(O_K)$.
  If $N$ is trivial, then $M$ is the standard log structure of $\Spec(O_K)$. 
  If $N$ is not trivial, then $M$ is not an fs log structure, but 
$M$ is always the filtered union of fs log structures contained in $M$ 
whose restriction on $\eta$ is $N$. 
  Below, the category of (weak) log abelian varieties for $M$ is defined to be the inductive limit of the categories of 
(weak) log abelian varieties for these fs log structures, and the polarizability for the log structure $M$ can be understood as that 
for some of these fs log structures. 
\end{para}

\begin{para} 
  The following categories were introduced in \cite{KKN4} 13.4. 

  Let $\cC_0^{\ptpol}$ be the category of log abelian varieties over $O_K$ with respect to the log structure $M$. 

  Let $\cC_1^{\ptpol}$ be  the category of pointwise polarizable log $1$-motifs over $O_K$ with respect to $M$.
 
  Let $\cC_0^{\pol}, \cC_1^{\pol}$ be the full subcategory of $\cC_0^{\ptpol}, \cC_1^{\ptpol}$ consisting of polarizable objects, respectively.

  Let $(\text{LAV}/K)$ be the category of log abelian varieties over $K$ with respect to the log structure $N$. 

  Let $\cC_2^{\ptpol}=\cC_2^{\pol}$ 
be the full subcategory of  $(\text{LAV}/K)$ consisting of objects  having the following property: if $[Y\to G_{\log}]$ denotes the corresponding log $1$-motif 
over $K$ and $0\to T\to G\to B\to 0$ denotes the exact sequence with $T$ a torus and $B$ an abelian variety, then $Y$ and $T$ are unramified and $B$ is of semistable reduction.
\end{para}

\begin{para}
We have the natural functors $\alpha_i\colon \cC_0^{\ptpol}\to \cC_i^{\ptpol}$ ($i=1,2$) 
(cf.\ \cite{KKN4} 13.2). 
  In particular, the functor $\alpha_1$ associates to a log abelian variety over $O_K$ 
the log 1-motif over $O_K$ defined by the family of the log $1$-motifs corresponding to the induced log abelian variety over 
$O_K/m_K^{n+1}$ for $n\geq 0$, where $m_K$ is the maximal ideal of $O_K$. 

  The functor $\alpha_1$ sends $\cC_0^{\pol}$ into $\cC_1^{\pol}$.
\end{para}

  We will prove the following results in Section \ref{s:dvr2}, which were announced in \cite{KKN4} 13.4. 

\begin{thm}
\label{t:cdvr}
$(1)$  $\cC_0^{\ptpol}=\cC_0^{\pol}$.

\medskip

$(2)$ The functor $\alpha_1$ induces an equivalence of categories $\cC_0^{\pol} \simeq \cC_1^{\pol}$.

\medskip

$(3)$ The functor $\alpha_2$ induces an equivalence of categories $\cC_0^{\pol} \simeq \cC_2^{\pol}$. 
\end{thm}

In particular, in the case $N$ is trivial, we have

\begin{cor} Taking the generic fiber, we have an equivalence from the category of log abelian varieties over $O_K$ for the standard log structure to the category of abelian varieties over $K$ with semistable reduction. 
\end{cor}

  First we have the following category equivalence. 
  The proof shows that a formal polarization induces a polarization at the generic fiber. 
  This implies a general fact used later in this paper that if a candidate of a polarization is a polarization at a point $s$ of the base, it is a polarization also at any generization $t$ of $s$ (cf.\ Lemma \ref{l:pol}). 

\begin{thm}\label{eqv1}  
We have an equivalence of categories
$$\cC_1^{\pol}\simeq \cC_2^{\pol}.$$
\end{thm}

  This is deduced from the corresponding result Theorem 15.10 in \cite{KKN4} by taking care of polarization. 
  In fact, for an object of $\cC_1^{\pol}$, it is easy to see that the corresponding log $1$-motif $[Y_2 \to G_{2,\log}]$ 
over $K$ defined in Section 15 in \cite{KKN4} is polarizable. 
  Thus we have a functor 
$$\gamma:\cC_1^{\pol}\to \cC_2^{\pol}.$$ 
  Conversely, let 
$p: [Y_2\to G_{2,\log}]\to [X_2\to G^*_{2,\log}]$ be a polarization on a log $1$-motif over $K$. 
Then $p$ induces a corresponding morphism $p: [Y_1\to G_{1,\log}]\to [X_1\to G^*_{1,\log}]$ over $O_K$ by Section 15 in \cite{KKN4}, which is easily seen 
to be a polarization. 
  Thus we have an inverse functor $$\delta: \cC_2^{\pol} \to \cC_1^{\pol}.$$ 

  The next is proved by Proposition 16.12 in \cite{KKN4} immediately. 

\begin{lem}\label{lem4} $(1)$ $\alpha_1\simeq \delta\alpha_2$ on $\cC_0^{\ptpol}$. 

\medskip

$(2)$ $\alpha_2\simeq \gamma\alpha_1$ on $\cC_0^{\pol}$.
\end{lem}

\section{GAGF for log abelian varieties, I}
\label{sec:GAGF1}

The aim of this section is to prove that a log abelian variety is determined by its formal completion  (Theorem \ref{thm5.1}). 

\begin{para}
  Let $(R,m)$ be a complete noetherian local ring endowed with an fs log structure on $\Spec(R)$. 

  Let $\cP$ be the category of polarizable $m$-adic formal log abelian varieties over $R$. 
  Here an $m$-adic formal object means a  family of objects $A_n$ over $R/m^{n+1}$ ($n\geq 0$) endowed with isomorphisms $A_{n}\otimes_{R/m^{n+1}} R/m^n\cong A_{n-1}$. 
  
  Note that the following proposition is easily seen. 
\end{para}

\begin{prop}
\label{prop8}
  We have an equivalence of categories
$$\cP'\overset{\cong}\to \cP,$$
where $\cP'$ is the category of polarizable log $1$-motifs over $R$.
The functor $\cP'\to \cP$ is given by taking the log abelian variety over $R/m^{n+1}$ corresponding to the induced polarizable log $1$-motif over $R/m^{n+1}$. 
\end{prop} 

  Let $\cQ$ be the category of log abelian varieties over $R$ and let $\cQ^{\pol}$ be the full subcategory of $\cQ$ consisting of polarizable objects.

\begin{thm}\label{thm5.1}
  Let $(R,m)$ be a complete noetherian local ring endowed with an fs log structure on $\Spec(R)$. 
  Then 
 the functor from $\cQ^{\pol}$ to $\cP$ 
is fully faithful. 
\end{thm}

\begin{para}
  In the rest of this section, we prove this theorem. 
  Let $A$ and $A'$ be objects of $\cQ^{\pol}$. 
  Assume that a morphism $A\to A'$ is given formally. 
  We have to prove that it is uniquely algebraized. 
  For simplicity, assume that it is an isomorphism.
  The general case is similar.
  Let $\Sig$ be the first standard fan (\cite{KKN5} 1.9) induced by a polarization of 
the pullback of $A$ to the closed point of $\Spec(R)$.
  This is a complete and wide fan by Proposition 10.3 of \cite{KKN5}. 
  Fix a prime number $\ell$ which is invertible on the base. 

  Let $n\ge 0$ be a nonnegative integer. 
  Let $(P_n,e_n,G_n,Q_n,\Sigma)$ be the object of $\mathcal B$ over $R/m^{n+1}$ corresponding to ($A_n := A\otimes R/m^{n+1}, \Sig)$ 
by Theorem \ref{t:A,B}. 
  This model $P_n$ is projective by Theorem 1.11 of \cite{KKN5}.
  On the other hand,
let $(P,e,G,Q,\Sigma)$ be the object of $\mathcal B$ corresponding
to $A$. 

  It suffices to show that
$(P_n,e_n,G_n,Q_n,\Sigma)$ $(n \ge0)$ determines $(P,e,G,Q,\Sigma)$. 
  For, by Theorem \ref{t:A,B}, 
it means that $A$ is recovered from $A_n$ $(n \ge0)$. 

  We show this. 
  Let $(P'_n,e'_n,G'_n,Q_n,\Sigma)$ and $(P',e',G',Q,\Sigma)$ be similar objects given by $A'$. 
  First, we are given a formal isomorphism $P^{\circ} \to (P')^{\circ}$ of formal schemes. 
  Here $(-)^{\circ}$ means to forget the log structure. 
  By GAGF of Grothendieck, we have an algebraic isomorphism $P^{\circ} \overset{\cong}\to (P')^{\circ}$ of schemes.
  Below we identify $P^{\circ}$ and $(P')^{\circ}$.
\end{para}

\begin{para}
  Next, we take care of log structures as follows.
  To this end, we prove $(\tilde P)^{\circ} \overset{\cong}\to (\tilde P')^{\circ}$. 

Consider the $\bar Y$-torsor on the \'etale site of $P^{\circ}=(P')^{\circ}$ given by local morphisms into $(\tilde P)^{\circ}$ (resp.\ $(\tilde P')^{\circ})$ over $P^{\circ}$. Since the restriction to the special fiber of these torsors are isomorphic, they are isomorphic by the proper base change theorem for $H^1$.
  Hence we can identify $(\tilde P')^{\circ}$ with $\tilde P^{\circ}$. 
\end{para}

\begin{para}
  We compare log structures of $P$ and that of $P'$. 
  We have a surjection $X\oplus M_S^{\gp}/\cO^\times_S\to M^{\gp}/\cO^\times$, where $M$ is the log structure of $\tilde P$.
  This is because the fan is constructed in $\Hom(\cS, \N)\times \Hom(X, \Z)$. 
  We have the following observations 1 and 2.

1. For each $x\in X$, we have a line bundle $L(x)$ which is the inverse image of the image of $x$ in $M^{\gp}/\cO^{\times}$ on $(\tilde P)^{\circ}$ and a line bundle $L'(x)$ on $(\tilde P')^{\circ}=(\tilde P)^{\circ}$. The actions of $\bar Y$ changes $L(x)$ and $L'(x)$, but $L(x)^{-1}L'(x)$ is unchanged and descends to $P^{\circ}$. By Grothendieck GAGF for line bundles,  the formal isomorphism $1\cong L(x)^{-1}L'(x)$ on $P^{\circ}$ becomes an algebraic isomorphism $1\cong L(x)^{-1}L'(x)$. 
Hence we have an isomorphism $L(x) \cong L'(x)$ on $(\tilde P)^{\circ}$ compatible with the action of $\bar Y$. 

2. We prove that $M/\cO^\times$ of $\tilde P$ coincides with that of $\tilde P'$. In fact, $M/\cO^\times$ of $\tilde P$ (resp.\ $\tilde P'$)  is a submonoid of a quotient group of $X\oplus M_S^{\gp}/\cO^\times_S$. The coincidence is checked as follows at each $t\in P$. Take a point $u$ of the special fiber of $\tilde P$ which belongs to the closure of $t$. Then since $\Spec(\hat \cO_{\tilde P,u})\to \Spec(\cO_{\tilde P, u})$ is surjective, the coincidence at $t$ can be checked by the coincidence at $\hat \cO_{P, u}$ and hence is checked formally.

These 1 and 2 prove that the two log structures $M$, $M'$ on $(\tilde P)^{\circ}$ have an isomorphism $M\cong M'$ which is compatible with the actions of $\bar Y$. Hence we have $M\cong M'$ on $P^{\circ}$. 
\end{para}

\begin{para}
  Thus $(P_n)$ $(n \ge0)$ determines a projective fs log scheme $P$.
  Again by GAGF for schemes, 
  $(e_n)$ $(n \ge0)$ determines a section $e$.

  Next, $P_n \to Q_n$ determines $P \to Q$. 
  Hence, $G$ is also recovered as the inverse image of the unit section 
of $Q$ by the last map. 

  The rest is the recovery of the action $G \times P \to P$ of $G$ on $P$ and the partial group law 
$(P \times P)' \to P$.
  The action is a part of the partial group law, and, since $(P_n \times P_n)'$ $(n \geq0)$ is 
represented by another projective scheme, again by GAGF, the partial 
group law is recovered. 
  This completes the proof of Theorem \ref{thm5.1}. 
\end{para}

\section{GAGF for log abelian varieties, II}\label{sec:GAGF2}

  Let the notation be as in Section \ref{sec:GAGF1}.

\begin{thm}\label{thm5} 
  Let $(R,m)$ be a complete noetherian local ring with an fs log structure on $\Spec(R)$. 
Then we have an equivalence of categories
$$\cQ^{\pol}\overset{\simeq}\to \cP.$$
\end{thm}

\begin{sbrem}
  The proof below also shows the GAGF for polarized objects, not only for polarizable objects. 
  Further, let $I$ be any ideal of $R$. 
  Then we can ask if we have the $I$-adic GAGF.
\end{sbrem}

\begin{para}
  In the rest of this section, we prove this theorem.
  Since we already show the full faithfulness in Theorem \ref{thm5.1}, it is enough to show that a given polarized formal log abelian variety can be algebraized. 

  First, we have an admissible pairing $X\times Y\to {\bold G}_{m,\log}/{\bold G}_m$, where $X$ and $Y$ are the ones for the closed fiber. 
  This is given $m$-adic formally, and extends automatically to $\Spec(R)$. Let $\bar X$ be the image of $X\to \cH om(Y, {\bold G}_{m,\log}/{\bold G}_m)$ and define $\bar Y$ similarly.

  Let $n\ge3$ be an integer invertible on the base. 
  To use Theorem \ref{t:localmoduli}, we give the abelian part an $n$-level structure. 
  More precisely, we algebraize the abelian part and after \'etale localizing the base if necessary, take an $n$-level structure on it. 
  By using the local moduli over $\bZ[1/n]$, which is log regular by Theorem \ref{t:localmoduli}, and by the full faithfulness (Theorem \ref{thm5.1}), 
we are reduced to the case where the base $S=\Spec(R)$ is log regular. 

  In the following, we assume that $S$ is log regular. 
  The outline is as follows. 
  We algebraize an object of $\cB$. 
  We then get the corresponding object of $\cA$ by Theorem \ref{t:A,B}.
  Finally we take care of polarizations. 
\end{para}

\begin{para}
\label{algebraizeP}
  As in the previous section, we take the first standard fan $\Sig$ for the closed fiber. 
  Then we have a formal object of $\cB$. 
  From this, we obtain a projective scheme $P$ by Grothendieck GAGF. 
  Since the base $S$ is log regular, the open set $U$ of $S$ where the log structure is trivial is dense.
  We endow $P$ with the log structure by the complement of the inverse image of $U$.
  Then we can check formally that $P$ is log smooth fs log scheme over $S$. 
  $Q$ extends to $S$, and $P\to Q$ is defined. 
  Further, the partial group law $(P \times P)' \to P$ is obtained again by GAGF of Grothendieck from the partial group laws $(P_n \times P_n)' \to P_n$ $(n \geq0)$. 
\end{para}

\begin{para}
  Let $G$ be the inverse image of the $0$-section by $P\to Q$.
  The group law of $G$ and the action of $G$ on $P$ are given by the partial group law of $P$. 
  In particular, we have a section $e$. 
\end{para}

\begin{para} 
  We prove that $P$ is a $G$-torsor, that is, we have $G\times P \cong P\times_Q P\; ;\; (a, x) \mapsto (x, ax)$. 

Let $(P\times P)''$ be the part of $P\times P$ consisting of $(a, b)$ such that $ab^{-1}$ in $Q$ belongs to $\Sigma$. 
  Then we have a morphism $(P\times P)''\to P\; ;\;(a,b)\mapsto a^{-1}b$.

  We have $$(P\times P)'\cong  (P\times P)''\; \quad (a,b) \mapsto (b, ab), \;(a^{-1}b, a)\leftarrow  (a,b).$$
This induces $G\times P \cong P\times_Q P$. 
\end{para}

\begin{para}
\label{semiabel}
  We prove that $G$ is semiabelian as in the following steps. 

1.  $G$ is representable, smooth and separated. 

  This is shown by taking the fan $\Sig'$ associated to a star. 
  By applying the argument in \ref{algebraizeP} to $\Sig \sqcap \Sig'$, we have another proper model $P'$.
  Since $G$ is the $\sig$-part of $P'$ so that it is an open of $P'$ and hence is representable and separated, where $\sig=\Hom(\cS, \bN) \times\{0\}$. 
  Further, $P'$ is log smooth and hence $G$ is log smooth. Since $G$ is strict over the base, $G$ is smooth. 

2.  $G$ is connected. 

  It is enough to show that $G$ is divisible.  
  By the construction in the proof of Theorem \ref{t:A,B}, we have an abelian sheaf $\tilde A$ and an exact sequence
$$0 \to G \to \tilde A \to \tilde Q \to 0.$$
  After the kfl sheafification, by \cite{Kato:FI2} Theorem 3.1, 
we have an exact sequence 
$$0 \to G \to (\tilde A)^{\kfl} \to (\tilde Q)^{\kfl} \to 0$$
of kfl sheaves.  
  (In actual, $\tilde A = (\tilde A)^{\kfl}$ but it is not necessary.) 
  Since $(\tilde A)^{\kfl}$ is divisible and $(\tilde Q)^{\kfl}$ is torsion-free, 
$G$ is divisible. 
  
3.  Let $\eta$ be any point of the base $S$. 
  We prove that $G$ has no additive part at $\eta$. 
  In the fiber at $\eta$, $G_{\overline \eta}$ has the torus part $T$, the abelian part $B$, and the additive part. 
  Let $t=\dim(T)$, $c=\dim(B)$, $a$ the dimension of the additive part, and let $t'$ be the rank of the stalk $\bar Y_{\overline \eta}$. 
  Let $d$ be the relative dimension of $G$ over $S$. 

  Again by the construction in the proof of Theorem \ref{t:A,B}, we have an abelian sheaf $A=\tilde A/\overline Y$ and an exact sequence
$$0 \to G \to A \to Q \to 0.$$
  From this, we have a complex
$$0 \to 
T_{\ell}(T)
\to T_{\ell}(A)_{\overline \eta(\mathrm{k\acute{e}t})} \overset f \to \overline Y_{\overline \eta}\otimes \bZ_{\ell} \to0$$
of $\bZ_{\ell}$-modules, which is exact except at the middle term, and 
the isomorphism 
$$\Ker(f)/T_{\ell}(T)\cong T_{\ell}(B).$$
  Here $T_{\ell}$ denotes the $\ell$-adic Tate module. 
  Further, we can apply \cite{KKN4} Proposition 18.1 to our $A$ and the results there still hold for $A$. 
  See a complementary explanation below in \ref{A[n]}. 
  In particular, $A[\ell^n]$ for any $n$ is represented by an fs log scheme over $S$ which is finite over $S$ and k\'et locally constant. 
  These show the equality $t+t'+2c=2d$. 
  On the other hand, we have trivially $a+t+c=d$. 
  Hence $2a+2t+2c=2d$. 
  Hence $2a+t=t'$. 
  We have $a=0$ if we can prove $t' \le t$.

  We prove $t' \le t$.
  Similarly in the proof of Theorem \ref{t:A,B}, the last statement of \cite{KKN2} Theorem 7.3 (1) implies that $\tilde A_{\overline \eta}$ is the pushout of 
$T_{\log}^{\prime(\overline Y_{\overline \eta})} \gets T' \to G_{\overline \eta}$ for some homomorphism $T' \to G_{\overline \eta}$, where $T'=\cHom(\overline X_{\overline \eta}, \Gmlog)$. 
  By the argument which follows in the proof of Theorem \ref{t:A,B}, we can see from the representablity of $\tilde P_{\overline \eta}$  that this homomorphism $T' \to G_{\overline \eta}$ is injective.
  Hence $t' \le t$. 
\end{para}

\begin{para}
\label{A[n]}
  In the above, we applied \cite{KKN4} Proposition 18.1 to our $A$.
  Here we remark on how to modify the proofs in \cite{KKN4}. 
  The difference between the current situation and that in \cite{KKN4} lies in that we only know that $A^{(\Sig)}$ is representable only for a specific $\Sig$. 

  For this reason, first, in the proof of \cite{KKN4} Lemma 18.4, which says that $A[n]$ ($n \ge1$) is represented by an algebraic space with an fs log structure, 
the representability of $A^{(n^{-1}\Sig_0)}$ is not trivial.
  Here $\Sig_0$ is the fan consisting of all translations of the cone $\sig$ in \ref{semiabel} 1.
  We prove that $A^{(n^{-1}\Sig_0)}$ is representable.
  In the construction till \ref{semiabel}, $\tilde A$, $A$, and $G$ do not depend on the choices of the complete and wide fan $\Sig$.
  In fact, for another complete and wide fan $\Sig'$, the fan $\Sig \sqcap \Sig'$ is also complete and wide, so that, to see this, we may assume that $\Sig'$ is a subdivision of $\Sig$.  Then, $\tilde P$ for $\Sig'$ coincides with the $\Sig'$-part of $\tilde A$ for $\Sig$. 
  From this, we see that $\tilde A$ for $\Sig'$ and that for $\Sig$ coincide, $A$'s coincide, and $G$'s also. 
  We return to the current situation with a fixed $\Sig$. 
  Since $n^{-1}\Sig$ is another complete and wide fan, $A^{(n^{-1}\Sig)}$ is representable.
  Then, $A^{(n^{-1}\Sig_0)}$ is representable as an open of a log blow-up of 
$A^{(n^{-1}\Sig)}$.
  Thus, from the argument in the proof of \cite{KKN4} Lemma 18.4, $A[n]$ is representable. 

  Next, we do not have the analogue of \cite{KKN4} Lemma 16.10 now, which was used to show that $A[n]$ is finite.
  But in the argument in \cite{KKN4} 16.11, where we use \cite{KKN4} Lemma 16.10, we need it only for one complete fan, and we can use our $\Sig$, for which the conclusion of \cite{KKN4} Lemma 16.10 holds.  
  Hence, we can also show that $A[n]$ is finite. 

  The remaining part of the proof is not necessarily to be changed. 
\end{para}

\begin{para}
  Thus we have obtained an object $(P, e, G, Q, \Sigma)$ of $\cB$. 
  By Theorem \ref{t:A,B}, this gives a weak log abelian variety $A$ with an exact sequence $0\to G\to A\to Q\to 0$. 
\end{para}

\begin{para} The rest is to show that the formal polarization becomes algebraic and that our $A$, which is a weak log abelian variety, is a log abelian variety and polarized. 
  Roughly speaking, the proof goes as follows. 
  On the model, the formal polarization gives a formal $\Gm$-torsor. This becomes algebraic by classical GAGF. 
  Take the associated $\Gmlog$-torsor. 
  Precisely, we argue as follows. 

  Let $\sig$ be the wide cone in $\Sig$. 
  Consider $I:=\sig$-part of $\tilde A$  and $J:=\ell^{-1}\sig$-part of $\tilde A$ (see the proof of Proposition 12.8 of \cite{KKN5}). 
  First, as in \cite{KKN5} Section 5, we have a formal $\Gm$-torsor on $P$ (pullback by the diagonal and take a special section of $\Gmlog/\Gm$-torsor on this special model) and algebraize it by the classical GAGF. 
Restrict it to $I$ and we have an algebraic $\Gm$-torsor so (by the extension of scalars) a $\Gmlog$-torsor $L$ on $I$. 
  By Proposition \ref{p:ffonLAV}, $L^{\otimes \ell^2}|_J$ descends to $J/G[\ell] \cong I$ and isomorphic to $L$. 
  Hence, by the argument in the proof of \cite{KKN5} Proposition 12.8, we have a $\Gmlog$-torsor on $\tilde A$, which descends to $A$. 

  Notice here that the associated biextension of this last torsor on $A$ formally coincides not with $p$ but with $2p$. 
  To do with this difference, we can proceed as follows. 
  Since $A$ satisfies the conditions in 1.4.1 of \cite{KKN5} (this is checked formally by \cite{KKN5} 4.14), by \cite{KKN5} Proposition 2.3, we have 
$$\Biext(A,A;\Gmlog) \cong \Hom(A,\cExt(A,\Gmlog)).$$ 
  We identify these two groups. 
  Let $$q\colon A \to \cExt(A,\Gmlog)$$ be the algebraization of $2p$. 
  Then, $q(A[2])=0$ because it is so formally. 
  (Here $[2]$ means the $2$-torsions.) 
  Consider the sequence of kfl sheaves $$0 \to A[2] \to A \overset 2 \to A \to 0.$$ 
  (By \cite{Z},  
any weak log abelian variety over a noetherian fs log scheme is a kfl sheaf.)
  Here the homomorphism $2\colon A \to A$ is surjective.  
  This is seen as follows. 
  Covering $A$ by copies of $\tilde A^{(\sig)}$ as in \ref{ketpre}, it is reduced to the fact that 
the morphism $2\colon \tilde A^{(2^{-1}\sig)} \to \tilde A^{(\sig)}$ is kfl surjective. 
  Thus the above sequence is exact. 
  Hence there is a homomorphism 
$$q'\colon A \to \cExt(A,\Gmlog)$$ such that $2q' = q$. 
  Since $q' - p$ is killed by 2 formally and $2$ is surjective, $q'=p$ formally.  
  Further, the biextension $q'$ is symmetric.
  This is checked formally, that is, checked in $\cH^1(A\times A, \Gmlog)$ by the use of Proposition \ref{p:ffonLAV}.
  Since $q'$ is a polarization at closed points, by Lemma \ref{l:pol} below, it is a polarization. 
  Therefore, we conclude that $q'$ is an algebraization of $p$. 
  This completes the proof of Theorem \ref{thm5}. 
\end{para}

\begin{lem}
\label{l:pol}
  Let $A$ be a weak log abelian variety over a noetherian fs log scheme $S$. 
If a symmetric biextension $p$ on $A$ is a polarization at a point $s\in S$, it is a polarization at any generization of $s$. 
\end{lem}

\begin{pf}
  By the proof of Theorem \ref{eqv1}, the lemma is valid if the underlying scheme of $S$ is the spectrum of a complete discrete valuation ring. 
  Then it holds in general by the reduction to this case.
\end{pf}

\section{Weak log abelian varieties over complete discrete valuation rings, II}
\label{s:dvr2}
  Here we prove Theorem \ref{t:cdvr}. 
  First, by varying fs log structures contained in $M$, we obtain from Theorem \ref{thm5} 
the following corollary, which proves (2) of Theorem \ref{t:cdvr}. 

\begin{cor}
\label{cor} 
  The functor $\alpha_1\colon \cC_0^{\pol}\simeq \cC_1^{\pol}$ gives an equivalence of categories.
\end{cor}

Together with Lemma \ref{lem4} (2), we have Theorem \ref{t:cdvr} (3). 

\begin{para}
  Finally, we prove Theorem \ref{t:cdvr} (1), that is, all objects of $\cC_0^{\ptpol}$ are polarizable.

  We denote by $\theta\colon \cC_1^{\pol}\simeq \cC_0^{\pol}$ the
  inverse functor of $\alpha_1\colon \cC_0^{\pol}\simeq \cC_1^{\pol}$
  in Corollary \ref{cor}, that is, an analogue of the Mumford construction. 

  Let $A$ be an object of $\cC_0^{\ptpol}$. 
  Then $\delta \alpha_2(A)$ is polarizable. Since $\alpha_1(A)$ is isomorphic to $\delta\alpha_2(A)$ by Lemma \ref{lem4} (1), $\alpha_1(A)$ is polarizable. It remains to prove $A\cong \theta\alpha_1(A)$. 
  By $\alpha_1\theta\simeq 1$ on $\cC_1^{\pol}$, we have 
$\alpha_1(A)\cong \alpha_1\theta \alpha_1(A)$. 
  Since $\alpha_1\colon \cC_0^{\ptpol}\to \cC_1^{\ptpol}$ is fully
  faithfull by Theorem \ref{thm5.1}, 
the last isomorphism implies $A\cong \theta\alpha_1(A)$. 
\end{para}

\noindent Takeshi Kajiwara

\noindent Department of Applied mathematics \\
Faculty of Engineering \\
Yokohama National University \\
Hodogaya-ku, Yokohama 240-8501 \\
Japan

\noindent kajiwara@ynu.ac.jp%
\par\bigskip\par
\noindent Kazuya Kato%

\noindent 
Department of Mathematics
\\
University of Chicago
\\
5734 S.\ University Avenue
\\
Chicago, Illinois, 60637 \\
USA
\\
\noindent kkato@math.uchicago.edu
\par\bigskip\par
\noindent Chikara Nakayama%

\noindent Department of Economics \\ Hitotsubashi University \\
2-1 Naka, Kunitachi, Tokyo 186-8601 \\ Japan

\noindent c.nakayama@r.hit-u.ac.jp

\end{document}